\documentclass[12pt]{article}
\usepackage{amsthm, amsmath, amssymb,verbatim,xspace, graphicx}
\usepackage{pdfsync, listings, color,pgf,tikz}
\usepackage[all,cmtip]{xy}
\usepackage{tikz}
\usepackage{breqn}
\usepackage{authblk}
\usepackage[T1]{fontenc}
\usepackage[utf8]{inputenc}
\usepackage{authblk}

\newcommand{\assignmentNumber}[1]{}

\newcommand\R{{\mathbb R}}

\providecommand{\keywords}[1]
{
  \small	
  \textbf{\textit{Keywords:}} #1
}

\def \beq{\begin{equation}}
	\def \eeq{\end{equation}}

\renewcommand{\rq}[1]{(\ref{#1})}

\theoremstyle{plain}

\title{Inverse Dynamic Problem for the Dirac System on Finite Metric Tree Graphs and the Leaf Peeling Method}
\author[1]{Sergei Avdonin\thanks{saavdonin@alaska.edu}}
\author[1]{Nina Avdonina\thanks{navdonina@alaska.edu}}
\author[2]{Olha Sus\thanks{olha.sus@tufts.edu}}
\affil[1]{Department of Mathematics and Statistics, University of Alaska Fairbanks}
\affil[2]{Department of Education, Tufts University}

\begin{document}
\maketitle

\begin{abstract} In this paper, we consider the inverse dynamic problem for the Dirac system on finite metric tree graphs. Our main goal is to recover the topology (connectivity) of a tree, lengths of edges, and a matrix potential function on each edge. We use the dynamic response operator as our inverse data and apply the \textit{leaf peeling method}. In addition, we present a new dynamic algorithm to solve the forward problem for the Dirac system on general finite metric graphs.
\end{abstract}


\section{Introduction}
\quad
 Under metric graphs or differential equation networks (DENs) we understand differential operators on geometric graphs coupled by certain vertex matching conditions (for more details see \cite{Kuchment}). Its theory plays a fundamental role in many sciences and engineering. The range for the applications of DENs is wide (see, for instance, \cite{2020_Avdonin_Edward}).

In recent years, the study of the Dirac operator on metric graphs has generated a growing interest (see, e.g., \cite{2016_Adamyan}, and the list of references therein). 
It has been used to build models for electronic properties of graphene \cite{2009_Neto_Peres}, propagation of electromagnetic waves in graphene-like photonic crystals \cite{2015_Dietz}, ultracold matter in optical lattices and many others (see, e.g., \cite{2021_Borrelli_Carlone_1}, \cite{2021_Borrelli_Carlone_2}, and references therein).

The Dirac operator on compact and non-compact graphs has been applied to spatially one-dimensional quantum mechanical problems (see, e.g., \cite{1990_Bulla_Trenkler}). It also has its application in solid state physics and optics (see \cite{2001_Rakesh} and the papers cited therein). 

There is a considerable amount of literature dedicated to solving inverse problems for the Schr\"{o}dinger equation on finite metric tree graphs (FMTG). For instance, in \cite{2008_Avdonin_Kurasov} three different inverse problems (spectral, dynamical and scattering) were studied and the \textit{leaf peeling method} (LPM) was proposed. This method was extended to solve boundary inverse problems for various types of PDE's on finite metric tree graphs and various matching conditions (see, e.g., \cite{2008_Avdonin_Kurasov}, \cite{2010_Avdonin_Kurasov_Nowaczyk}, \cite{2010_Avdonin_Leugering}, \cite{2015_Avdonin_Bell}, \cite{2020_Avdonin_Zhao}). There are two versions of the LP method. The spectral version is when the unknown data is recovered from the spectral inverse data (the Titchmarsh-Weyl matrix function), while the dynamical is when the unknown data is recovered from the dynamical inverse data (dynamic response operator). The former approach is less consuming in terms of mathematical computations, whereas the latter one is more beneficial from the programming point of view since it contains only integral equations and lower triangular matrices. 

Let us remark that both versions of LPM, the spectral and the dynamic, are recursive. It allows one to recalculate the spectral or dynamical data for the peeled tree from the original tree. It is realized by peeling the leaves (leaf edges) from the original tree to form a new peeled tree for which the corresponding reduced response data is recalculated. 

While the inverse problems for various types of PDE's on FMTG have been actively studied, there are quite a few papers dedicated to solving the inverse problem for 1-D Dirac systems on FMTG. In paper \cite{2014_Belishev_Mikhailov}, the inverse problem for a one-dimensional dynamical Dirac system on a half-line was solved, where the authors used an approach which required the extension of the BC-method to the time-domain inverse problem for the Dirac system. In paper \cite{2018_Mikhaylov_Murzabekova}, the authors considered the inverse dynamic and spectral problems for the 1-D Dirac system on a finite tree. The aim of the paper was to recover the topology of a tree (lengths and connectivity of edges) as well as the matrix potentials on each edge by applying a spectral version of LPM.

In this paper, we solve the inverse dynamic problem for the 1-D Dirac system on FMTG by applying the new, comprehensive version of the LP method proposed in \cite{2020_Avdonin_Zhao}. It allows us to recover the topology (connectivity) of a tree, lengths of edges, and a matrix potential function on each edge. We do not use any spectral methods at any step of our solution. To our best knowledge this is the first paper solving the inverse problem for the 1-D Dirac system on finite metric tree graphs by applying a pure dynamic algorithm. 

In addition, we propose a new dynamic algorithm for solving the forward problem for the 1-D Dirac system on general finite metric graphs. The original version of this algorithm for solving the initial boundary value problem for the 1-D wave equation on general graphs was introduced in \cite{2019_Avdonin_Zhao}. It was further applied to solve the exact controllability problem for the wave equation on graphs in \cite{2022_Avdonin_Zhao}.

The outline of the paper is as follows. In Section 2, we introduce the definition of the finite metric tree graph $\Omega$ and the 1-D Dirac system associated with $\Omega$. In Section 3, we prove the existence of the unique generalized solution to the forward problem for the 1-D Dirac system on general metric graphs. In Section 4, we solve the inverse problem for the 1-D Dirac system on finite metric tree graphs by recovering all unknown data. 

\section{Preliminaries}
\quad
Let $\Omega=(V,E)$ be a finite metric graph. By $E$ we denote the set of edges $\{e_1,\cdots,e_N\}$ connected at the vertices $V=\{v_1,\cdots,v_M\}$. Every edge $e_k\in E$ is identified with an interval $(0,l_k)$ of the real line. The notation $e_k\sim v$ means that the edge $e_k$ is incident to the vertex $v$. By $\Gamma=\{\gamma_1,\cdots,\gamma_L\}$ we denote the boundary of $\Omega$, that is the set of vertices having multiplicity one, and we assume that $\gamma_L$ is a root of $\Omega$. In addition, on each edge $e_k$, we identify the vertex which is located further from the root vertex $\gamma_L$ as $x=0$, and the closer vertex to the root as $x=l_k$.

Let us also introduce the following notations. Let $J:=\begin{pmatrix}0&1\\-1&0\end{pmatrix}$ be the Pauli matrix. Let $Q=\begin{pmatrix}p_k&q_k\\q_k&-p_k\end{pmatrix}$ be a real matrix-valued potential with $p_k, q_k\in C^1(e_k;\mathbb{R})$ defined on each edge $e_k$, $k=1,\cdots,N$. We denote by $H=L^2(\Omega;\mathbb{C}^2):=\bigoplus_{k=1}^N L^2(e_k;\mathbb{C}^2)$ a complex valued space of square integrable functions on the graph $\Omega$, endowed with the standard norm $
\|\cdot\|^2_{L^2(\Omega;\mathbb{C}^2)}:=\sum_{k=1}^N\|\cdot\|^2_{L^2(e_k;\mathbb{C}^2)}$.
For the element $U\in L^2(\Omega;\mathbb{C}^2)$ we write $U:=(u^1,u^2)^{\top}\big|_{e_k}=\{(u_k^1,u_k^2)^{\top}\}_{k=1}^N$, $u_k^1,u_k^2\in L^2(e_k;\mathbb{C})$.

We now associate the following initial boundary value problem (IBVP) for the 1--D Dirac system with the graph $\Omega$:

\begin{equation}\label{Equation1}
iU_t+JU_x+QU=0,\quad t\geq0,\quad x\in e_k,\quad k=1,\cdots,N,
\end{equation}
\begin{equation}\label{Equation2}
u_k^1(v,t)=u_s^1(v,t),\quad e_k\sim v,\quad e_s\sim v,\quad k\neq s, \quad v\in V\backslash \Gamma,\quad t\geq0,
\end{equation}
\begin{equation}\label{Equation3}
\sum_{k|e_k\sim v}u_k^{2,\pm}(v,t)=0, \quad v\in V\backslash \Gamma,\quad t\geq0,
\end{equation}
\begin{equation}\label{Equation4}
u^1(\gamma_k,t)=f_k(t),\quad k=1,\cdots,L-1,\quad u^1(\gamma_L,t)=0,\quad t\in[0,T],
\end{equation}
\begin{equation}\label{Equation5}
U(\cdot,0)=0,
\end{equation}
where $T$ is arbitrary positive number, $F=(f_1(t),\cdots,f_{L-1}(t))^{\top}\in L^2(0,T;\mathbb{R}^{L-1})$ is a boundary control vector function. Also, here $u_k^{2,\pm}(v,t)$ stands for $u_k^2(0,t)$ or $-u_k^2(l_k,t)$ according to whether $x$ is equal to $0$ or $l_k$ at vertex $v$, $k=1,\cdots,N$.

In the next section we present a new dynamic approach for solving the forward problem for the 1-D Dirac system on general metric graphs.

\section{Forward problem for the Dirac system on general metric graphs}

\quad
First, we present the solution to the forward problem for the 1-D Dirac system on a finite interval and a star-shaped graph.

\subsection{Forward problem for the Dirac system on a finite interval}

Let $Q=\begin{pmatrix}p&q\\q&-p\end{pmatrix}$ be a real matrix-valued potential defined on a finite interval $(0,l)$. 

We consider the following IBVP with a control function $f(t)$ applied at the left-hand end of $(0,l)$:
\begin{equation}\label{Equation8}
iU_t+JU_x+QU=0,\quad 0<x<l,\quad 0<t<T,
\end{equation}
\begin{equation}\label{Equation9}
u^1(0,t)=f(t),\quad u^1(l,t)=0,\quad 0<t<T,
\end{equation}
\begin{equation}\label{Equation10}
U(x,0)=0,\quad 0<x<l.
\end{equation}

We use the notation $U^{f,+}(x,t)=\begin{pmatrix}u^{1,f,+}(x,t)\\u^{2,f,+}(x,t)\end{pmatrix}$ to represent the solution to system (\ref{Equation8})-(\ref{Equation10}).

In paper \cite{2014_Belishev_Mikhailov}, the corresponding forward problem for the 1-D Dirac system was considered  on the semi-axis $x>0$. The authors showed that the defined on semi-axis matrix-valued potential $Q$ uniquely determines  the vector-kernel function $(w^{1,+},w^{2,+})^{\top}$ such that, for any $T>0,$ the solution to the 1-D Dirac system on the semi-axis has the following Duhamel-type representation:
\begin{equation}\label{Equation11}
U^{f,+}(x,t)=\begin{cases}
\begin{pmatrix}1\\i\end{pmatrix}f(t-x)+\int\limits_x^t\begin{pmatrix}w^{1,+}(x,s)\\w^{2,+}(x,s)\end{pmatrix}f(t-s)\,ds,\quad x\geq 0,\quad 0\leq t\leq T,\\
0,\quad 0<t<x,\end{cases}
\end{equation}
where $(w^{1,+},w^{2,+})^{\top}|_{t<x}=0$, $(w^{1,+},w^{2,+})^{\top}|_{\Delta^T}\in C^1(\Delta^T;\mathbb{C}^2)$, $w^{1,+}(0,\cdot)=0$ $\bigl(\Delta^T:=\{(x,t)~\big|~x\geq 0,~0\leq t\leq T,~x\leq t\}\bigr)$.

Using this representation we can  find the solution of the system (\ref{Equation8})-(\ref{Equation10}) for all $T>0.$
 We extend  potential functions $p$ and $q$ to the semi-axis $x>0$ by the following rule:
\begin{equation*}
p(2nl\pm x)=p(x),\quad q(2nl\pm x)=q(x),\quad \text{for all}\,\,n\in\mathbb{N},
\end{equation*}
and hence the extension of the matrix-valued potential $Q$ to semi-axis satisfies the equalities
\begin{equation} \label{qex}
Q(2nl\pm x)=Q(x),\quad \text{for all}\,\,n\in\mathbb{N}.
\end{equation}

The extended matrix-valued potential $Q$ uniquely determines the vector-kernel function $(w^{1,+},w^{2,+})^{\top}$ for $x>0$. Therefore, by taking into account (\ref{Equation11}) and \rq{qex}
we state the following proposition which can be proved by direct substitution. 

\vspace{0.1cm}

\textbf{Proposition 3.1.1.} If $p,q\in C^1_{loc}([0,\infty);\mathbb{R})$, $f\in L^2([0,T];\mathbb{R})$, and $t\geq0$, then the IBVP (\ref{Equation8})-(\ref{Equation10}) has a unique generalized solution $U^{f,+}\in C([0,T];L^2([0,T];\mathbb{C}^2))$ given by the formula
\begin{equation}\label{Equation12}
\begin{split}
U^{f,+}(x,t)=\sum_{n=0}^{\lfloor\frac{t-x}{2l}\rfloor}\left(\begin{pmatrix}1\\i\end{pmatrix}f(t-2nl-x)+\int\limits_{2nl+x}^t\begin{pmatrix}w^{1,+}(2nl+x,s)\\w^{2,+}(2nl+x,s)\end{pmatrix}f(t-s)\,ds\right)\\
+\sum_{n=1}^{\lfloor\frac{t+x}{2l}\rfloor}\left(\begin{pmatrix}-1\\i\end{pmatrix}f(t-2nl+x)+\int\limits_{2nl-x}^t\begin{pmatrix}-w^{1,+}(2nl-x,s)\\w^{2,+}(2nl-x,s)\end{pmatrix}f(t-s)\,ds\right),
\end{split}
\end{equation}
where $\lfloor\cdot\rfloor$ is a floor-function. 
Here and in the rest of the paper we assume that $f$ is extended to the whole real line and $f(t)=0$ for $t<0$.

{For the scalar wave equation similar representation of the solution was obtained in \cite{2019_Avdonin_Zhao}.

Now we consider the IBVP on a finite interval $(0,l)$ with a control function $f$ applied at the right-hand end:
\begin{equation}\label{Equation13}
iU_t+JU_x+Q(x)U=0,\quad 0<x<l,\quad t>0,
\end{equation}
\begin{equation}\label{Equation14}
u^1(0,t)=0,\quad u^1(l,t)=f(t),\quad t>0,
\end{equation}
\begin{equation}\label{Equation15}
U(x,0)=0,\quad 0<x<l.
\end{equation}
We use the notation $U^{f,-}(x,t)=\begin{pmatrix}u^{1,f,-}(x,t)\\u^{2,f,-}(x,t)\end{pmatrix}$ to represent the solution to the system (\ref{Equation13})-(\ref{Equation15}).

To solve this IBVP  we construct $\widetilde{Q}(x)=Q(l-x)$ with $\tilde{p}(x)=p(l-x)$, $\tilde{q}(x)=q(l-x)$ and extend $\widetilde{Q}(x)$ to the semi-axis by letting
\begin{equation*}
\widetilde{Q}(2nl\pm x)=\widetilde{Q}(x),\quad \text{for all}\,\,n\in\mathbb{N}.
\end{equation*}

The extended matrix-valued potential $\widetilde{Q}$ uniquely determines the vector-kernel function $(w^{1,-},w^{2,-})^{\top}$ for $x>0$. The analog of the Proposition 3.1.1 is
valid for $U^{f,-},$ and the following representation can be checked by direct substitution:
\begin{equation}\label{Equation16}
\begin{aligned}
U^{f,-}(x,t)={} & \sum_{n=0}^{\lfloor\frac{t+x-l}{2l}\rfloor}\Biggl(\begin{pmatrix}1\\-i\end{pmatrix}f(t-(2n+1)l+x)\\
 & +\int\limits_{(2n+1)l-x}^t\begin{pmatrix}w^{1,-}((2n+1)l-x,s)\\-w^{2,-}((2n+1)l-x,s)\end{pmatrix}f(t-s)\,ds\Biggr)\\
 & -\sum_{n=0}^{\lfloor\frac{t-x-l}{2l}\rfloor}\Biggl(\begin{pmatrix}1\\i\end{pmatrix}f(t-(2nl+1)l-x)\\
 & +\int\limits_{(2n+1)l+x}^t\begin{pmatrix}w^{1,-}((2n+1)l+x,s)\\w^{2,-}((2n+1)l+x,s)\end{pmatrix}f(t-s)\,ds\Biggr).
\end{aligned}
\end{equation}

\subsection{Forward problem for the Dirac system on a star-shaped graph}
\quad
In this subsection we consider the forward problem which is associated with the 1-D Dirac system on a star-shaped graph $\Omega$. Here $\Omega=\{V,E\}$ with $V=\{v,\gamma_1,\cdots,\gamma_L\}$. Each edge $e_k$, $k=1,\cdots,L$, is identified with the interval $(0,l_k)$, the vertex $v$ is identified with $x=0$ and $\gamma_k$ is identified with $x=l_k$.

We now associate the following IBVP for the 1–D Dirac system with the star-shaped graph $\Omega$: 
\begin{equation}\label{star_1}
iU_t+JU_x+QU=0,\quad t\geq 0,\quad x\in e_k,\quad k=1,\cdots,L,
\end{equation}
\begin{equation}\label{star_2}
u_k^1(0,t)=u_s^1(0,t),\quad k\neq s,\quad t\geq0,\quad k,s=1,\cdots,L,
\end{equation}
\begin{equation}\label{star_3}
\sum_{k=1}^Lu_k^{2,\pm}(0,t)=0,\quad t\geq0,
\end{equation}
\begin{equation}\label{star_4}
u_k^1(l_k,t)=f_k(t),\quad k=1,\cdots,L-1,\quad u_L^1(l_L,t)=0,\quad 0\leq t\leq T,
\end{equation}
\begin{equation}\label{star_5}
U(\cdot,0)=0.
\end{equation}

Let us denote the common value $u_k^1(x,t)$ at $x=0$ by $g(t)$ for $k=1,\cdots,L$. We have that $g(t)=0$ for $t<l_1$.

\begin{figure}[ht]
\begin{center}
\begin{tikzpicture}
%% vertices
\draw[fill=black] (0,0) circle (1.5pt);
\draw[fill=black] (4,0) circle (1.5pt);
\draw[fill=black] (2,1) circle (1.5pt);
\draw[fill=black] (2,3) circle (1.5pt);
\draw[fill=black] (4,2) circle (1.5pt);
%% vertex labels
\node at (-0.5,0) {$l_1,\gamma_1$};
\node at (4.6,0) {$l_2,\gamma_2$};
\node at (1.2,1.2) {$v, g(t)$};
\node at (2,0.7) {$0$};
\node at (2,3.3) {$l_L,\gamma_L$};
\node at (4,2.4) {$l_{L-1},\gamma_{L-1}$};
%%% edges
\draw[thick] (0,0) -- (2,1) -- (2,3) -- (2,1) -- (4,0);
\draw[dashed] (4,2) -- (2,1);
\end{tikzpicture}
\caption{Star-shaped graph $\Omega$}
\label{star-shaped-graph}
\end{center}
\end{figure}
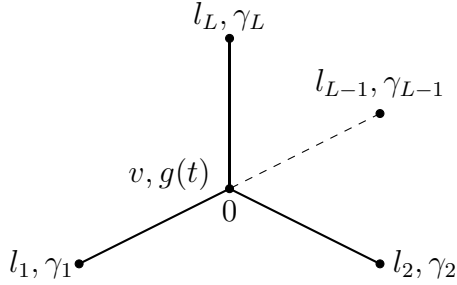

By superposition principle, one obtains the following formulas for the solution to the problem (\ref{star_1})-(\ref{star_5}):
\begin{equation}\label{Equation17}
\begin{split}
U_k(x,t)=U_k^{f_k,-}(x,t)+U_k^{g,+}(x,t),\quad k=1,\cdots,L-1,\\
U_L(x,t)=U_L^{g,+}(x,t).
\end{split}
\end{equation}

To find an explicit representation of the solution $U_k(x,t)$ for each $k=1,\cdots,L$, one needs to solve for the unknown function $g(t)$. Using formulas (\ref{Equation16}) and (\ref{Equation12}), for each $k=1,\cdots,L$ one obtains:
\begin{equation*}
\begin{aligned}
U_k(x,t)={}
& \sum_{n=0}^{\lfloor\frac{t+x-l_k}{2l_k}\rfloor}\Biggl(\begin{pmatrix}1\\-i\end{pmatrix}f_k(t-2nl_k-l_k+x)
+\int\limits_{2nl_k+l_k-x}^t\begin{pmatrix}w_k^{1,-}(2nl_k+l_k-x,s)\\-w_k^{2,-}(2nl_k+l_k-x,s)\end{pmatrix}f_k(t-s)\,ds\Biggr)\\
& -\sum_{n=0}^{\lfloor\frac{t-x-l_k}{2l_k}\rfloor}\Biggl(\begin{pmatrix}1\\i\end{pmatrix}f_k(t-2nl_k-l_k-x) +\int\limits_{2nl_k+l_k+x}^t\begin{pmatrix}w_k^{1,-}(2nl_k+l_k+x,s)\\w_k^{2,-}(2nl_k+l_k+x,s)\end{pmatrix}f_k(t-s)\,ds\Biggr)\\
& +\sum_{n=0}^{\lfloor\frac{t-x}{2l_k}\rfloor}\Biggl(\begin{pmatrix}1\\i\end{pmatrix}g(t-2nl_k-x)
+\int\limits_{2nl_k+x}^t\begin{pmatrix}w_k^{1,+}(2nl_k+x,s)\\w_k^{2,+}(2nl_k+x,s)\end{pmatrix}g(t-s)\,ds\Biggr)\\
& +\sum_{n=1}^{\lfloor\frac{t+x}{2l_k}\rfloor}\Biggl(\begin{pmatrix}-1\\i\end{pmatrix}g(t-2nl_k+x)
+\int\limits_{2nl_k-x}^t\begin{pmatrix}-w_k^{1,+}(2nl_k-x,s)\\w_k^{2,+}(2nl_k-x,s)\end{pmatrix}g(t-s)\,ds\Biggr), 
\end{aligned}
\end{equation*}
where $k=1,\cdots,L-1$, and for $k=L$ we have
\begin{equation*}
\begin{split}
U_L(x,t)=\sum_{n=0}^{\lfloor\frac{t-x}{2l_L}\rfloor}\left(\begin{pmatrix}1\\i\end{pmatrix}g(t-2nl_L-x)+\int\limits_{2nl_L+x}^t\begin{pmatrix}w_L^{1,+}(2nl_L+x,s)\\w_L^{2,+}(2nl_L+x,s)\end{pmatrix}g(t-s)\,ds\right)\\
+\sum_{n=1}^{\lfloor\frac{t+x}{2l_L}\rfloor}\left(\begin{pmatrix}-1\\i\end{pmatrix}g(t-2nl_L+x)+\int\limits_{2nl_L-x}^t\begin{pmatrix}-w_L^{1,+}(2nl_L-x,s)\\w_L^{2,+}(2nl_L-x,s)\end{pmatrix}g(t-s)\,ds\right).
\end{split}
\end{equation*}

Now, substituting these expressions to the Kirchhoff-type conditions (\ref{star_2})-(\ref{star_3}), we obtain the equation for the unknown function $g(t)$:
\begin{equation*}
\begin{split}
Lig(t)+\int\limits_0^t\left(\sum_{k=1}^Lw_k^{2,+}(0,s)\right)g(t-s)\,ds\\
+2\sum_{k=1}^{L-1}\sum_{n=0}^{\lfloor\frac{t-l_k}{2l_k}\rfloor}\Bigl(-if_k(t-2nl_k-l_k)-
\int\limits_{2nl_k+l_k}^tw_k^{2,-}(2nl_k+l_k,s)f_k(t-s)\,ds\Bigr)\\+
2\sum_{k=1}^{L-1}\sum_{n=1}^{\lfloor\frac{t}{2l_k}\rfloor}\Bigl(ig(t-2nl_k)+
\int\limits_{2nl_k}^tw_k^{2,+}(2nl_k,s)g(t-s)\,ds\Bigr)\\
+
2\sum_{n=1}^{\lfloor\frac{t}{2l_L}\rfloor}\Bigl(ig(t-2nl_L)+
\int\limits_{2nl_L}^tw_L^{2,+}(2nl_L,s)g(t-s)\,ds\Bigr)=0.
\end{split}
\end{equation*}

From this it follows that
\begin{equation}\label{Equation18}
Lig(t)+\int\limits_0^t\left(\sum_{k=1}^Lw_k^{2,+}(0,s)\right)g(t-s)\,ds=G(t),
\end{equation}

where
\begin{equation*}
\begin{split}
G(t)=2\sum_{k=1}^{L-1}\sum_{n=0}^{\lfloor\frac{t-l_k}{2l_k}\rfloor}\Bigl(if_k(t-2nl_k-l_k)+
\int\limits_{2nl_k+l_k}^tw_k^{2,-}(2nl_k+l_k,s)f_k(t-s)\,ds\Bigr)\\
-2\sum_{k=1}^L\sum_{n=1}^{\lfloor\frac{t}{2l_k}\rfloor}\Bigl(ig(t-2nl_k)+
\int\limits_{2nl_k}^tw_k^{2,+}(2nl_k,s)g(t-s)\,ds\Bigr).
\end{split}
\end{equation*}

Equation (\ref{Equation18}) is a delay integral equation with respect to the unknown function $g(t)$. We will prove now that this equation can be solved by method of steps.

\vspace{0.1cm}

\textbf{Proposition 3.2.1.} Suppose that $f_k(t)\in L^2_{loc}([0,\infty);\mathbb{R})$ 
for $k=1,\cdots,L-1$ in the IBVP (\ref{star_1})-(\ref{star_5}). Then one can compute $g(t)$ from (\ref{Equation18}), and $g(t)\in L^2_{loc}([0,\infty);\mathbb{R})$.  

\begin{proof}
Let $l_{min}=\min\limits_{k=1,\cdots,L}l_k$. We have that the RHS of (\ref{Equation18}) depends on $g(s)$ and $f_k(s)$, $k=1,\cdots,L-1$, where $g(s)$ is delayed by at least $t-2l_{min}$ and $f_k(s)$ is known for all $s$ values. From this it follows that if $g$ is known for $s\in[0,t-2l_{min}]$, then for $s\in[t-2l_{min},t]$ we obtain a Volterra integral equation of the second kind with respect to $g(t)$. Since we know that $g(t)=0$ for $t\leq l_{min}$, one can solve (\ref{Equation18}) iteratively with a time step size of $2l_{min}$.
\end{proof}
Since $g(t)$ is known, then for each $k=1,\cdots,L$ one can obtain an explicit representation of a solution $U_k(x,t)$ to the stated problem (\ref{star_1})-(\ref{star_5}) on a star-shaped graph $\Omega$.

In the next subsection we will solve forward problem associated with the 1-D Dirac system on general graphs.

\subsection{Forward problem for the Dirac system on general graphs}

\quad
Let $\Omega=(V,E)$ be a finite metric graph associated with (\ref{Equation1})-(\ref{Equation5}), where $V=\{v_1,\cdots,v_M\}$ is the set of vertices of $\Omega$ and $E=\{e_1,\cdots,e_N\}$ is the set of edges. Let us denote by $g_i(t)=u^1(v_i,t)$ for each $v_i\in V$ and let $g(t)=(g_1(t),\cdots,g_{M}(t))^T$. We denote the set of boundary vertices by $\Gamma=\{\gamma_1,\cdots,\gamma_L\}$ and $g(t)$ on $\Gamma$ is defined as $g(t)=u^1(\gamma_i,t)$, $i=1,\cdots,L$. The values of $g(t )$ on $\Gamma$ are given as boundary conditions in (\ref{Equation4}). Thus it remains to compute $g(t)$ on $V\backslash \Gamma$.

Let $e_k=(v_i,v_j)$ be an edge identified with the interval $[0,l_k]$, where $k=1,\cdots,N$. We suppose that $v_i$ is identified with the left-hand end of the interval $[0,l_k]$. Let us define the operators $W^{\pm}_k:L^2(0,T)\mapsto C([0,T];L^2(0,l_k))$ by the following rule:
\begin{equation*}
\begin{split}
(W_k^{+}f)(x,t)=U_k^{f,v_i}(x,t)=(u_k^{1,f,v_i},u_k^{2,f,v_i})^{\top},\\
(W_k^{-}f)(x,t)=U_k^{f,v_j}(x,t)=(u_k^{1,f,v_j},u_k^{2,f,v_j})^{\top},
\end{split}
\end{equation*}
where $U_k^{f,v_i}(x,t)$ represents the solution on the edge $e_k$ with the control function applied at $x=v_i$ and $U_k^{f,v_j}(x,t)$ represents the solution on the edge $e_k$ with the control function applied at $x=v_j$.

Let us define an $N\times M$ matrix operator $A$ such that $\{a_{ki}\}_{k,i=1}^{N,I}=W_k^{+}$, $\{a_{kj}\}_{k,j=1}^{N,J}=W_k^{-}$ $(I+J=M)$ if there is an edge $e_k$ from $v_i$ to $v_j$ and $\{a_{km}\}_{k,m=1}^{N,M}=0$ otherwise. As one can observe, $Ag$ gives us a column vector, where the $k$-th entry equals to $U_k(x,t)$ on the edge $e_k$. If one knows the vector function $g(t)$, then $Ag$ is the solution to the stated forward problem (\ref{Equation1})-(\ref{Equation5}).

Now let us define observation operators $\mathcal{O}_k^{+}$ and $\mathcal{O}_k^{-}$ on $C([0,T];L^2(0,l_k))$ by the following rule:
\begin{equation*}
\mathcal{O}_k^{+}U_k=u_k^2(v_i,\cdot),
\end{equation*}
\begin{equation*}
\mathcal{O}_k^{-}U_k=-u_k^2(v_j,\cdot).
\end{equation*}

Then, we obtain four combinations of $\mathcal{O}^{\pm}_k$ and $W^{\pm}_k$ on $e_k$ for $k=1,\cdots,N$:
\begin{equation}\label{composition1}
\begin{split}
(\mathcal{O}_k^{-}W_k^{+})(f)=\mathcal{O}^{-}_{k}U_k^{f,v_i}(x,t)=-u_k^{2,v_i}(v_j,t)\\
=-\sum_{n=0}^{\lfloor\frac{t-l_k}{2l_k}\rfloor}\left(if(t-2nl_k-l_k)+\int\limits_{2nl_k+l_k}^tw_k^{2,+}(2nl_k+l_k,s)f(t-s)\,ds\right)\\
-\sum_{n=1}^{\lfloor\frac{t+l_k}{2l_k}\rfloor}\left(if(t-2nl_k+l_k)+\int\limits_{2nl_k-l_k}^tw^{2,+}_k(2nl_k-l_k,s)f(t-s)\,ds\right).
\end{split}
\end{equation}

\begin{equation}\label{composition2}
\begin{split}
(\mathcal{O}_k^{+}W_k^{+})(f)=\mathcal{O}^{+}_{k}U_k^{f,v_i}(x,t)=u_k^{2,v_i}(v_i,t)\\
=if(t)+\int\limits_0^tw_k^{2,+}(0,s)f(t-s)\,ds\\
+2\sum_{n=1}^{\lfloor\frac{t}{2l_k}\rfloor}\left(if(t-2nl_k)+\int\limits_{2nl_k}^tw^{2,+}_k(2nl_k,s)f(t-s)\,ds\right).
\end{split}
\end{equation}

\begin{equation}\label{composition3}
\begin{split}
(\mathcal{O}_k^{-}W_k^{-})(f)=\mathcal{O}^{-}_{k}U_k^{f,v_j}(x,t)=-u_k^{2,v_j}(v_j,t)\\
=\sum_{n=0}^{\lfloor\frac{t}{2l_k}\rfloor}\left(if(t-2nl_k)+\int\limits_{2nl_k}^tw_k^{2,-}(2nl_k,s)f(t-s)\,ds\right)\\
+\sum_{n=0}^{\lfloor\frac{t-2l_k}{2l_k}\rfloor}\left(if(t-2nl_k-2l_k)-\int\limits_{2nl_k+2l_k}^tw^{2,-}_k(2nl_k+2l_k,s)f(t-s)\,ds\right).
\end{split}
\end{equation}

\begin{equation}\label{composition4}
\begin{split}
(\mathcal{O}_k^{+}W_k^{-})(f)=\mathcal{O}^{+}_{k}U_k^{f,v_i}(x,t)=u_k^{2,v_j}(v_i,t)\\
=-2\sum_{n=0}^{\lfloor\frac{t-l_k}{2l_k}\rfloor}\left(-if(t-2nl_k-l_k)-\int\limits_{2nl_k+l_k}^tw^{2,-}_k(2nl_k+l_k,s)f(t-s)\,ds\right).
\end{split}
\end{equation}

By $B$ we denote the matrix operator of size $M\times N$ such that $\{b_{ik}\}_{i,k=1}^{I, N}=\mathcal{O}_k^{-}$, $\{b_{jk}\}_{j,k=1}^{J,N}=\mathcal{O}_k^{+}$ $(I+J=M)$ if there is an edge $e_k$ from $v_i$ to $v_j$ and $\{b_{km}\}_{k,m=1}^{N,M}=0$ otherwise. From this it follows that $BAg$ is a column vector of $M$ entries and the $i$-th entry represents the sum of second components of a vector function $U_k$ at a vertex $v_i$ (over all incident edges taken outwards of $v_i$).

Since the Kirchhoff-type conditions (\ref{Equation2})-(\ref{Equation3}) hold only at internal vertices $v_i$, we will define a diagonal $M\times M$ matrix $D$ to pick out the interior vertices by the following rule:
\begin{equation*}
\{d_{kj}\}_{k,j=1}^M=1,\quad\text{if}~k=j,~v_k\in V\backslash \Gamma\quad\text{and}\quad\{d_{kj}\}_{k,j=1}^M=0\quad\text{otherwise}.
\end{equation*}

Hence, taking into account the Kirchhoff-type conditions (\ref{Equation2})-(\ref{Equation3}), we obtain
\begin{equation}\label{Forward_KN}
DBAg=0.
\end{equation}

The above equation (\ref{Forward_KN}) is a system of $|V\backslash\Gamma|$ equations similar to (\ref{Equation18}). Indeed, let us introduce the following notations: 
$$
I(v_m):=\{e_j\in E\,|\,e_j\sim v_m,\,\,j=1,\cdots,N\},\quad m=1,\cdots,|V\backslash\Gamma|,
$$
and
$$
J(\gamma_i):=\{e_i\in E\,|\,e_i\sim \gamma_i,\,\,i=1,\cdots,L-1\}.
$$

We then obtain the number of $|V\backslash\Gamma|$ Volterra integral equations of the second kind of the following form
\begin{equation}\label{Forward_Volterra}
i\deg{(v_m)}g_m(t)+\int\limits_0^t\mathcal{K}_m(0,s)g_m(t-s)\,ds=G_m(t),
\end{equation}
where
$$
\mathcal{K}_m(0,s)=\sum_{k\in I(v_m)}w^{2,+}_k(0,s)
$$
are known and
\begin{equation*}
\begin{split}
G_m(t)=-\sum_{k\in J(v_m)\cap J(v_j)}\Bigl[\sum_{n=0}^{\lfloor\frac{t-l_k}{2l_k}\rfloor}\Bigl(if(t-2nl_k-l_k)+\int\limits_{2nl_k+l_k}^tw_k^{2,+}(2nl_k+l_k,s)f(t-s)\,ds\Bigr)\\
+\sum_{n=1}^{\lfloor\frac{t+l_k}{2l_k}\rfloor}\Bigl(if(t-2nl_k+l_k)+\int\limits_{2nl_k-l_k}^tw^{2,+}_k(2nl_k-l_k,s)f(t-s)\,ds\Bigr)\Bigr]\\
-2\sum_{k\in I(v_m)}\sum_{n=1}^{\lfloor\frac{t}{2l_k}\rfloor}\Bigl(ig_m(t-2nl_k)+
\int\limits_{2nl_k}^tw_k^{2,+}(2nl_k,s)g_m(t-s)\,ds\Bigr),
\end{split}
\end{equation*}
which contains a function $g_m(t)$ with a delay by at least $2\min\limits_{1\leq m\leq |V\backslash\Gamma|}l_m$.

The following proposition summarizes results obtained above.

\textbf{Proposition 3.3.1.} 
Suppose that $F\in L^2_{loc}([0,\infty);
\mathbb{R}^{L-1})$ in the IBVP (\ref{Equation1})-(\ref{Equation5}). Then, for $m=1,\cdots,|V\backslash \Gamma|$, one can compute $g_m(t)$ from (\ref{Forward_Volterra}), and $g_m(t)\in L^2_{loc}([0,\infty);\mathbb{R})$.  

\begin{proof}
See Proposition 3.2.1.
\end{proof}

Taking into account Proposition 3.3.1 and by applying the superposition principle we obtain the solution of the system (\ref{Equation1})-(\ref{Equation5}) on each edge $e_k$, $k=1,\cdots,N$ of the general metric graph $\Omega$ which is given by:
\begin{equation*}
U_k(x,t)=U_k^{g_i,+}(x,t)+U_k^{g_j,-}(x,t),
\end{equation*}
where $i,j\in\{1,\cdots,M\}$ and $i\neq j$.

In the next section we study the inverse problem for the 1-D Dirac system on finite metric tree graphs.

\section{Inverse problem for the Dirac system on finite metric tree graphs. LP method}
\quad
Let $\Omega=(V,E)$ be a finite metric tree graph with $V=\{v_1,\cdots,v_M\}$, $E=\{e_1,\cdots,e_N\}$, and $\Gamma=\{\gamma_1,\cdots,\gamma_L\}$, where $\Gamma$ is the set of boundary vertices of $\Omega$ with a root vertex $\gamma_L$.

In this section we solve the inverse dynamic problem for the 1-D Dirac system (\ref{Equation1})-(\ref{Equation5}) on $\Omega$. To realize it, we apply the LP algorithm.

\subsection{Dynamical and spectral inverse problems statements}

\quad
First, let us formulate the statement of the dynamical inverse problem (DIP) for the 1-D Dirac system (\ref{Equation1})-(\ref{Equation5}) on $\Omega$.

We introduce the response operator for the problem (\ref{Equation1})-(\ref{Equation5}) by the rule:
\begin{equation}\label{response_operator}
R^{T}\{F\}(t):=u^2(\cdot,t)\big|_{\Gamma\backslash\{\gamma_L\}}, \quad t\in[0,T],
\end{equation}
where $\Bigl(R^T\bigl(f^1(t),\cdots,f^{L-1}(t)\bigr)^{\top}\Bigr)(t)=\bigl(u^2(v_1,t),\cdots, u^2(v_{L-1},t)\bigr)^{\top}$.

The response operator $R^T$ has a form of convolution:
\begin{equation*}
(R^TF)(t)=(\mathbf{R}\ast F)(t)=\int\limits_0^t\mathbf{R}(t-s)F(s)\,ds,
\end{equation*}
where $\mathbf{R}(t)=\{R_{ij}\}_{i,j=1}^{L-1}$ is a response matrix. 

The entries of $\mathbf{R}(t)$ are defined in the following way: let $U_i$ be a solution to the IBVP (\ref{Equation1})-(\ref{Equation5}) with special boundary condition (\ref{Equation4}). That is,
\begin{equation*}
F=(0,\cdots,\underbrace{\delta(t)}_{\textrm{i-th}},\cdots,0)^{\top}.
\end{equation*}

Then
\begin{equation}\label{response_matrix}
R_{ij}(t)=u^2_i(v_j,t).
\end{equation}

Now, we can formulate the statement of the DIP associated with (\ref{Equation1})-(\ref{Equation5}).

\textbf{DIP:} Recover the tree $\Omega$ (connectivity of edges and their lengths) and the matrix potential function $Q$ on all edges of $\Omega$ from the response operator $R^T(t)$, given by (\ref{response_operator}), for $t>0$.

It is well-known (see, for instance, \cite{2008_Avdonin_Kurasov}, \cite{2017_Avdonin_Mikhaylov_Nurtazina}) that one can connect the dynamic problem (DP) (\ref{Equation1})-(\ref{Equation5}) to the spectral problem (SP) associated with $\Omega$. 

Indeed, let us put
$$
\Phi:=\begin{pmatrix}\varphi^1\\\varphi^2\end{pmatrix}\in L^2(\Omega),\quad \varphi^1_k, \varphi^2_k\in H^1(e_k),\quad k=1,\cdots, N
$$
and introduce the operator 
$$
\mathcal{L}\Phi:=J\Phi_x+Q\Phi,\quad x\in e_k
$$
with the domain
$$
dom(\mathcal{L})=\{\Phi\in L^2(\Omega):\,\,\varphi_k^1,\varphi_k^2\in H^1(e_k),\,\Phi\,\text{satisfies (\ref{Equation2}), (\ref{Equation3})},\,\varphi^1(v)=0,\,v\in\Gamma\}.
$$

The spectral problem (SP) associated with $\Omega$ has the following representation:
\begin{equation}\label{sp_1}
J\Phi_x+Q\Phi=\lambda\Phi,\quad x\in e_k, \quad k=1,\cdots,N,
\end{equation}
\begin{equation}\label{sp_2}
\varphi_k^1(v)=\varphi_s^1(v),\quad e_k\sim v,\quad e_s\sim v,\quad v\in V\backslash\Gamma,\quad k\neq s,
\end{equation}
\begin{equation}\label{sp_3}
\sum_{k|_{e_k\sim v}}\varphi_k^{2,\pm}(v)=0,\quad v\in V\backslash\Gamma,
\end{equation}
\begin{equation}\label{sp_4}
\varphi^1(v)=0,\quad v\in\Gamma.
\end{equation}

Now, for $F\in L^2(0,T;\mathbb{R}^{L-1})\cap(C^{\infty}_{0}(0,+\infty))^{L-1}$ we define the Fourier transform
\begin{equation}\label{F_transform}
\widehat{F}(k):=\int\limits_0^{\infty}F(t)e^{ikt}\,dt.
\end{equation}

Going formally in DP over the Fourier transform (\ref{F_transform}), we obtain SP with $\lambda=k$. 

Accordingly, one can obtain the connection between the spectral and dynamic inverse data (see, \cite{2008_Avdonin_Kurasov}, \cite{2017_Avdonin_Mikhaylov_Nurtazina}). We introduce the Titchmarsh-Weyl (T-W) matrix function in the following way: for $\lambda\not\in \mathbb{R}$ and $\xi\in \mathbb{R}^{L-1}$ we consider the problem (\ref{sp_1})-(\ref{sp_3}) with the non-homogeneous boundary condition:
\begin{equation}\label{sp_5}
\varphi^1(v_i)=\xi_i,\quad i=L-1.
\end{equation}

The TW matrix function connects the values of the solution $\Phi(\cdot,\lambda)$ to (\ref{sp_1})-(\ref{sp_3}), (\ref{sp_5}) in the way:
\begin{equation}\label{connection_1}
\varphi^2(\cdot,\lambda)|_{\Gamma}=\mathbf{M}(\lambda)\varphi^1(\cdot,\lambda)|_{\Gamma},
\end{equation}
\begin{equation}\label{connection_2}
\Bigl(\varphi^2(v_1,\lambda),\cdots,\varphi^2(v_{L-1},\lambda)\Bigr)^{\top}=\mathbf{M}(\lambda)(\xi_1,\cdots,\xi_{L-1})^{\top}.
\end{equation}

Thus, by applying the Fourier transform (\ref{F_transform}), we obtain
\begin{equation}\label{connection}
\mathbf{M}(k)=\int\limits_0^{\infty}\mathbf{R}(t)e^{ikt}\,dt,
\end{equation}
where (\ref{connection}) is understood in a weak sense.

Therefore, one can formulate the statement of the spectral inverse problem (SIP) associated with (\ref{sp_1})-(\ref{sp_4}).

\textbf{SIP:} Recover the tree (connectivity of edges and their lengths) and the matrix potential function $Q$ on all edges from $\mathbf{M}(\lambda)$.

\subsection{Auxiliary definitions and facts}

\quad
We need the following definitions and facts in order to solve the DIP.

\textbf{Definition 4.2.1.} In a star-shaped graph $\Omega^{\ast}$ all but one of its edges are called \textit{leaf} edges. The one edge that is not a leaf edge is called the \textit{stem} edge.

\vspace{0.1cm}

\textbf{Definition 4.2.2.} Let $S$ be a star-shaped subgraph of a tree graph $\Omega$. We call $S$ a \textit{sheaf} if it contains all edges of the graph incident to some internal vertex $v\in V\backslash\Gamma$ and if all but one of its edges are leaf edges (see Figure 2 below).
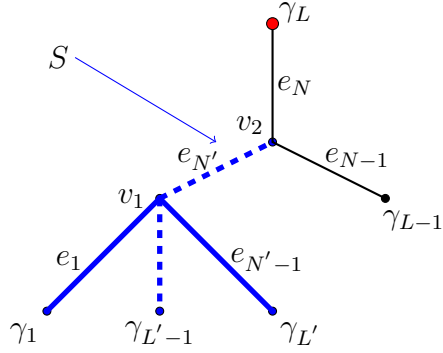
\begin{figure}[ht]
\begin{center}
\begin{tikzpicture}[scale=0.75]
%% vertices
\draw[fill=blue] (0,0) circle (2pt);
\draw[fill=black] (4,0) circle (2pt);
\draw[fill=blue] (2,1) circle (2pt);
\draw[fill=red] (2,3.1) circle (3pt);
\draw[fill=blue] (-2,-2) circle (2pt);
\draw[fill=blue] (2,-2) circle (2pt);
\draw[fill=blue] (0,-2) circle (2pt);
%% vertex labels
\node at (-0.5,0) {$v_1$};
\node at (-1.6,-1.1) {$e_1$};
\node at (1.9,-1.1) {$e_{N^{'}-1}$};
\node at (0.7,0.7) {$e_{N^{'}}$};
\node at (3.5,0.7) {$e_{N-1}$};
\node at (2.4,2) {$e_N$};
\node at (4.5,-0.4) {$\gamma_{L-1}$};
\node at (1.6,1.3) {$v_2$};
\node at (2.4,3.3) {$\gamma_L$};
\node at (-2.4,-2.4) {$\gamma_1$};
\node at (2.5,-2.4) {$\gamma_{L^{'}}$};
\node at (0,-2.4) {$\gamma_{L^{'}-1}$};
%%% edges
\draw[thick,blue,line width=2pt] (-2,-2)-- (0,0) -- (2,-2);
\draw[thick,blue,line width=2pt, dashed] (0,0) -- (0,-2);
\draw[thick] (2,3) -- (2,1) -- (4,0);
\draw[thick,blue,line width=2pt, dashed] (0,0) -- (2,1);

\node at (-1.8,2.5) {$S$};
\draw[->, blue, bend right] (-1.5,2.5) -- (1,1);

\end{tikzpicture}
\end{center}
\caption{Graph $\Omega$ and a sheaf $S$.}\label{Figure:figure_2}
\end{figure}

Additionally, we use the following proposition which shows the existence of a sheaf $S$ in a tree graph $\Omega$.

\vspace{0.1cm}

\textbf{Proposition 4.2.1} (see \cite{2015_Avdonin_Choque}). A tree graph, which is not a star-shaped graph, contains at least one sheaf.

\subsection{Recovering the tree and the potential matrix function on the leaves of a sheaf}

\quad
Let $U^{\delta}$ be a solution to the BVP for the 1-D Dirac system associated with  the tree graph $\Omega$: $U^{\delta}$ satisfies (\ref{Equation1}) on $e_k$, $k=1,\cdots,N$, and (\ref{Equation2}), (\ref{Equation3}) at the internal vertices $v\in V\backslash\Gamma$. At the boundary vertices $\gamma_s\in\Gamma$, $s=1,\cdots,L$, it satisfies the following conditions
\begin{equation*}
u^1(\gamma_1,t)=\delta(t),\quad u^1(\gamma_2,t)=0,\cdots, u^1(\gamma_L,t)=0.
\end{equation*}

Let us consider a sheaf $S$ in the tree $\Omega$ and assume that the boundary edge $e_1$ belongs to $S$. We denote by $l_1$ the length of the edge $e_1=[\gamma_1,v_1]$ which is identified with the finite interval $[0,l_1]$ and by $N^{'}$ the number of all edges in a sheaf $S$. We also assume that the edge $e_1$ is connected at the internal vertex $v_1\in V\backslash \Gamma$ to other $N^{'}-1$ edges $e_2,\cdots,e_{N^{'}}$ which are identified with the intervals $[l_1,l_1+l_k]$, where $l_k$ is the length of $e_k$, $k=2,\cdots,N^{'}$, $N^{'}\leq N$.

Now, let us denote by $g(t)=g^{v_1}:=u_k^1(v_1,t)$, $v_1\in V\backslash \Gamma$. For $t<3l_1$, from (\ref{Forward_Volterra}) and (\ref{Equation16}) we obtain
\begin{equation*}
\begin{split}
N^{'}ig(t)+\int\limits_0^t\left(\sum_{k=1}^{N^{'}}w_k^{2,+}(0,s)\right)g(t-s)\,ds\\
=-2i\delta(t-l_1)-2\int\limits_{l_1}^tw_1^{2,-}(l_1,s)\delta(t-s)\,ds=-2i\delta(t-l_1)+\psi_1(t),
\end{split}
\end{equation*}
where $\psi_1(t)$ is a piecewise continuous function.

Thus, we get
\begin{equation*}
ig(t)=-\frac{2}{N^{'}}i\delta(t-l_1)+\psi_2(t).
\end{equation*}

By taking into account the definition of the response operator, we obtain
\begin{equation*}
\begin{split}
R_{11}(t)=u_1^{\delta,2}(0,t)=u_1^{\delta,2,+}(0,t)+u_1^{g,2,-}(0,t)\\
=ig(t-l_1)+\int\limits_{l_1}^tw_1^{2,-}(l_1,s)g(t-s)\,ds+
i\delta(t)+\int\limits_0^tw_1^{2,+}(0,s)\delta(t-s)\\
+i\delta(t-2l_1)+\int\limits_{2l_1}^tw_1^{2,+}(2l_1,s)\delta(t-s)\,ds.
\end{split}
\end{equation*}

From this it follows that
\begin{equation*}
R_{11}(t)=i\delta(t)+i\frac{N^{'}-2}{N^{'}}\delta(t-2l_1)+\psi_3(t),\quad t<3l_1,
\end{equation*}
where $\psi_3(t)$ is some piecewise continuous function.

Hence, by knowing the entry $R_{11}(t)$, one can determine the length $l_1$ of the edge $e_1$ and the number of edges connected to $e_1$. To generalize, one can determine the lengths $l_k$ of the boundary edges $e_k$, $k=1,\cdots,N^{'}-1$ from the diagonal response entries $R_{ii}(t)$ in a sheaf $S$, where $i=1,\cdots,L^{'}$, $L^{'}\leq L-1$ and $L^{'}$ is the number of the boundary vertices in a sheaf $S$.

In paper \cite{2018_Mikhaylov_Murzabekova}, it was proved that the diagonal elements of the response matrix function $\mathbf{R}(t)$ determine not only the lengths of the boundary edges, but also the matrix potential function $Q$ defined on those edges. From this it follows that for a fixed boundary vertex $\gamma_j$, $j=1,\cdots,L^{'}$, $L^{'}\leq L-1$, one can recover potential functions $p_k$, $q_k$ on $e_k$ from $R_{kk}(t)$ for $t=2l_k$, $k=1,\cdots,N^{'}-1$, where $l_k$ is a length of $e_k$ which was recovered from $R_{kk}(t)$ previously.

Similarly, by identifying other sheaves in the tree graph $\Omega$, one can recover the rest of the lengths of the leaf edges of $\Omega$, the potential matrix function $Q$ determined on them, and the degree of the internal vertices which belong to those sheaves. 

\subsection{Spectral version of the LP method. Recovering the T-W function for a peeled tree}

\quad
Let $S$ be a sheaf in the tree graph $\Omega$. We assume that $S$ consists of boundary vertices $\{\gamma_1,\cdots,\gamma_{L^{'}}\}\in\Gamma$, the corresponding boundary edges $e_1,\cdots,e_{N^{'}-1}$ and an internal edge $e_{N^{'}}$. At this point we assume that we already recovered the lengths $l_k$ and potential functions $p_k$, $q_k$, $k=1,\cdots,N^{'}-1$ on boundary edges of $S$. We identify each edge $e_k$, $k=1,\cdots,N^{'}$ with the interval $[0,l_k]$ and the vertex $v_1$, the internal vertex of a sheaf, with the set of common endpoints $x=0$. For convenience, we renumerate the edge $e_{N^{'}}$ as $e_0$ and the vertex $v_1$ as $\gamma_0$.

By $\widetilde{\mathbf{M}}(\lambda)$ we denote the reduced T-W matrix function associated with the new peeled graph (see Figure 3 below) $\widetilde{\Omega}=\Omega\backslash\bigcup_{k=1}^{N^{'}-1}\{e_k\}$ with boundary points $\Gamma=\{\gamma_0,\gamma_1,\cdots,\gamma_{L-1}\}$.

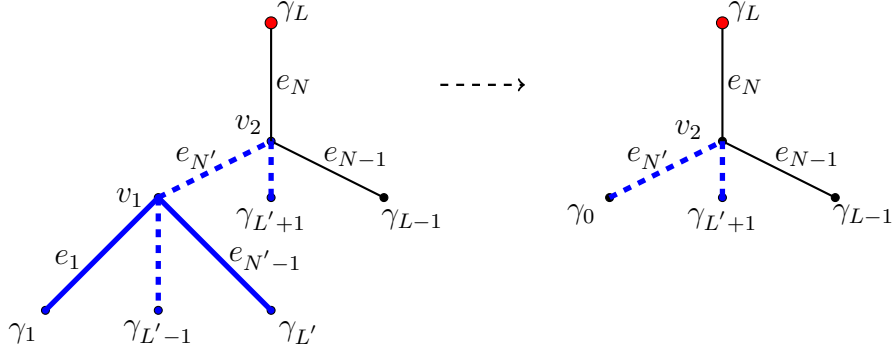
\begin{figure}[ht]
\begin{center}
\begin{tikzpicture}[scale=0.75]
%% vertices
\draw[fill=blue] (0,0) circle (2pt);
\draw[fill=black] (4,0) circle (2pt);
\draw[fill=blue] (2,1) circle (2pt);
\draw[fill=red] (2,3.1) circle (3pt);
\draw[fill=blue] (-2,-2) circle (2pt);
\draw[fill=blue] (2,-2) circle (2pt);
\draw[fill=blue] (0,-2) circle (2pt);
\draw[fill=blue] (2,0) circle (2pt);

%% vertex labels
\node at (-0.5,0) {$v_1$};
\node at (-1.6,-1.1) {$e_1$};
\node at (1.9,-1.1) {$e_{N^{'}-1}$};
\node at (0.7,0.7) {$e_{N^{'}}$};
\node at (3.5,0.7) {$e_{N-1}$};
\node at (2.4,2) {$e_N$};
\node at (4.5,-0.4) {$\gamma_{L-1}$};
\node at (1.6,1.3) {$v_2$};
\node at (2.4,3.3) {$\gamma_L$};
\node at (-2.4,-2.4) {$\gamma_1$};
\node at (2.5,-2.4) {$\gamma_{L^{'}}$};
\node at (0,-2.4) {$\gamma_{L^{'}-1}$};
\node at (2,-0.4) {$\gamma_{L^{'}+1}$};
%%% edges
\draw[thick,blue,line width=2pt] (-2,-2)-- (0,0) -- (2,-2);
\draw[thick,blue,line width=2pt, dashed] (0,0) -- (0,-2);
\draw[thick] (2,3) -- (2,1) -- (4,0);
\draw[thick,blue,line width=2pt, dashed] (0,0) -- (2,1);

\draw[->,line width=0.3mm,dashed] (5,2) -- (6.5,2);
\draw[thick,blue,line width=2pt, dashed] (2,1) -- (2,0);

\draw[fill=black] (8,0) circle (2pt);
\draw[fill=black] (12,0) circle (2pt);
\draw[fill=black] (10,1) circle (2pt);
\draw[fill=red] (10,3.1) circle (3pt);
\draw[fill=blue] (10,0) circle (2pt);
%% vertex labels
\node at (12.5,-0.3) {$\gamma_{L-1}$};
\node at (9.4,1.2) {$v_2$};
\node at (7.5,-0.3) {$\gamma_0$};
\node at (10.4,3.3) {$\gamma_L$};
\node at (8.7,0.7) {$e_{N^{'}}$};
\node at (11.5,0.7) {$e_{N-1}$};
\node at (10.4,2) {$e_N$};
\node at (10,-0.4) {$\gamma_{L^{'}+1}$};
%%% edges
\draw[thick] (10,3) -- (10,1) -- (12,0);
\draw[thick,blue,dashed,line width=2pt] (8,0) -- (10,1);
\draw[thick,blue,dashed,line width=2pt] (10,1) -- (10,0);

\end{tikzpicture}
\end{center}
\caption{Graph $\Omega$ and a peeled graph $\widetilde{\Omega}$.}\label{Figure:figure_3}
\end{figure}

First, we recalculate entries $\widetilde{M}_{0k}(\lambda)$, $k=0,L^{'}+1,\cdots,L-1$. Let us fix $\gamma_1$ which is the boundary point of the sheaf $S$. Let $\Phi$ be a solution to the problem (\ref{sp_1})-(\ref{sp_3}) with the following boundary conditions
\begin{equation*}
\varphi^1(\gamma_1)=1,\quad \varphi^1(\gamma_2)=0,\quad\cdots,\quad \varphi^1(\gamma_L)=0.
\end{equation*}

The function $\Phi$ solves the IVP on edge $e_1$:
\begin{equation}\label{ip_1}
J\Phi_x+Q\Phi=\lambda\Phi,\quad x\in e_1,
\end{equation}
\begin{equation}\label{ip_2}
\varphi^1(\gamma_1)=1,\quad \varphi^2(\gamma_1)=M_{11}(\lambda)
\end{equation}
and on edges $e_2,\cdots,e_{N^{'}-1}$ it solves the problem
\begin{equation}\label{ip_3}
J\Phi_x+Q\Phi=\lambda\Phi,\quad x\in e_k,
\end{equation}
\begin{equation}\label{ip_4}
\varphi^1(\gamma_k)=0,\quad \varphi^2(\gamma_k)=M_{1k}(\lambda),\quad k=2,\cdots,N^{'}-1.
\end{equation}

Since we know potentials functions $p_k$, $q_k$ on edges $e_1,\cdots,e_{N^{'}-1}$, we can solve problems (\ref{ip_1})-(\ref{ip_2}) and (\ref{ip_3})-(\ref{ip_4}), and use the K.-N. conditions (\ref{Equation2}), (\ref{Equation3}) at the internal vertex $\gamma_0$ of a sheaf $S$ to recover $\varphi_0^1(\gamma_0,\lambda)$ and $\varphi_0^2(\gamma_0,\lambda)$. Here $\Phi_0(\gamma_0,\lambda)=\bigl(\varphi_0^1(\gamma_0,\lambda),\varphi^2_0(\gamma_0,\lambda)\bigr)^{\top}$ is the value of the solution $\Phi$ at a new boundary point of a peeled tree $\widetilde{\Omega}$. Then we obtain:
\begin{equation*}
\widetilde{M}_{00}=\frac{\varphi_0^2(\gamma_0,\lambda)}{\varphi_0^1(\gamma_0,\lambda)},
\end{equation*}
\begin{equation*}
\widetilde{M}_{0k}(\lambda)=\frac{M_{1k}(\lambda)}{\varphi_0^1(\gamma_0,\lambda)},\quad k=L^{'}+1,\cdots,L-1. 
\end{equation*}

To find $\widetilde{M}_{k0}(\lambda)$, $k=L^{'}+1,\cdots,L-1$, we fix boundary point $\gamma_k$, $k\not\in\{1,\cdots,L^{'},L\}$ and consider the solution $\Phi$ to (\ref{sp_1})-(\ref{sp_3}) with the following boundary conditions
\begin{equation*}
\varphi^1(\gamma_k)=1,\quad \varphi^1(\gamma_s)=0,\quad k\neq s.
\end{equation*}

The function $\Phi$ solves the IVP on edges $e_1,\cdots,e_{N^{'}-1}$:
\begin{equation}\label{ip_5}
J\Phi_x+Q\Phi=\lambda\Phi,\quad x\in e_k,
\end{equation}
\begin{equation}\label{ip_6}
\varphi^1(\gamma_s)=0,\quad \varphi^2(\gamma_s)=M_{ks},\quad s=1,\cdots,N^{'}-1.
\end{equation}

Since we know the potential functions $p_k$, $q_k$ on edges $e_1,\cdots,e_{N^{'}-1}$, we can solve (\ref{ip_5})-(\ref{ip_6}) and use K.-N. conditions at the internal vertex $\gamma_0$ to recover values $\varphi_0^1(\gamma_0,\lambda)$ and $\varphi^2_0(\gamma_0,\lambda)$.

On the other hand, on the peeled graph $\widetilde{\Omega}$ the function $U$ solves the problem (\ref{sp_1})-(\ref{sp_3}) with the boundary conditions:
\begin{equation*}
\varphi^1(\gamma_k)=1,\quad \varphi^1(\gamma_0)=\varphi^1_0(\gamma_0,\lambda),
\end{equation*}
\begin{equation*}
\varphi^1(\gamma_s)=0,\quad s=L^{'}+1,\cdots, L,\quad \gamma_k\neq \gamma_s,\quad \gamma_s\neq \gamma_0.
\end{equation*}

Thus we obtain the relations:
\begin{equation*}
\widetilde{M}_{k0}=\varphi^2_0(\gamma_0,\lambda)-\varphi_0^1(\gamma_0,\lambda)\widetilde{M}_{00}(\lambda),
\end{equation*}
\begin{equation*}
\widetilde{M}_{ks}=M_{ks}-\varphi^1_0(\gamma_0,\lambda)\widetilde{M}_{0s}(\lambda).
\end{equation*}

To recover all elements of the reduced T-W matrix function, we need to repeat this procedure for all $k,s=L^{'}+1,\cdots,L-1$. 

Thus, using the described algorithm we can recalculate re reduced T-W matrix function for the new peeled tree $\widetilde{\Omega}$, that is, the tree without the leaf edges of a sheaf $S$.

Repeating this algorithm a sufficient number of times, we recover the tree $\Omega$ (connectivity of edges and their lengths) and the matrix potential function $Q$ defined on it.

\subsection{Dynamical version of the LP method. Recovering the response operator for a peeled tree}

\quad
We consider the dynamical problem (\ref{Equation1})-(\ref{Equation5}) on $\Omega$. Let $S$ be a sheaf in the tree graph $\Omega$. We assume that $S$ consists of boundary vertices $\{\gamma_1,\cdots,\gamma_{L^{'}}\}\in\Gamma$, the corresponding boundary edges $e_1,\cdots,e_{N^{'}-1}$ and an internal edge $e_{N^{'}}$. At this point we assume that we already recovered the lengths $l_k$ and potential functions $p_k$, $q_k$, $k=1,\cdots,N^{'}-1$ on boundary edges of $S$. We identify each edge $e_k$, $k=1,\cdots,N^{'}$ with the interval $[0,l_k]$ and the vertex $v_1$, the internal vertex of a sheaf, with the set of common endpoints $x=0$. For convenience, we renumerate the edge $e_{N^{'}}$ as $e_0$ and the vertex $v_1$ as $\gamma_0$.

We denote by $\widetilde{\mathbf{R}}(t)$ the reduced response function associated with the new peeled graph $\widetilde{\Omega}=\Omega\backslash\bigcup_{k=1}^{N'-1}\{e_k\}$ with boundary points $\Gamma=\{\gamma_0,\gamma_1,\cdots,\gamma_{L-1}\}$.

First, we recalculate entries $\widetilde{R}_{0k}(\lambda)$, $k=0,L^{'}+1,\cdots,L-1$. Let us fix $\gamma_1$ which is the boundary point of the sheaf $S$. We consider the function $U^{\delta}$ that is a solution of the problem (\ref{Equation1})-(\ref{Equation3}), (\ref{Equation5}) with the following boundary conditions
$$
u^{1,\delta}(\gamma_1,t)=\delta(t),\quad u^{1,\delta}(\gamma_k,t)=0,\quad k=2,\cdots, L.
$$

The function $U^{\delta}$ solves the following Cauchy problem on edge $e_1$:
\begin{equation}\label{dp-1}
\begin{cases}
iU^{\delta}_t+JU_x^{\delta}+QU^{\delta}=0,\quad x\in(0,l_1),\\
u^{1,\delta}(\gamma_1,t)=\delta(t),\quad t\in\mathbb{R},\\
u^{2,\delta}(\gamma_1,t)=R_{11}(t),\quad t>0,\\
U^{\delta}(x,0)=0,\quad x\in(0,l_1).
\end{cases}
\end{equation}

On other boundary edges of $S$, the function $U^{\delta}$ solves the problems:
\begin{equation}\label{dp-2}
\begin{cases}
iU^{\delta}_t+JU_x^{\delta}+QU^{\delta}=0,\quad x\in(0,l_k),\\
u^{1,\delta}(\gamma_k,t)=0,\quad t\in\mathbb{R},\\
u^{2,\delta}(\gamma_1,t)=\begin{cases}0,\quad t<l_1+l_k,\\R_{1k}(t),\quad t\geq l_1+l_k,\end{cases}\\
U^{\delta}(x,0)=0,\quad x\in(0,l_k),\quad k=2,\cdots,L^{'}.
\end{cases}
\end{equation}

Since potential functions $p_k$, $q_k$ are known on edges $e_1,\cdots,e_{N^{'}-1}$, we can solve problems (\ref{dp-1}) and (\ref{dp-2}). By using compatibility conditions (\ref{Equation2}) and (\ref{Equation3}) at the internal vertex $\gamma_0$, we can recover the values of $u_0^{1,\delta}(\gamma_0,t)$ and $u_0^{2,\delta}(\gamma_0,t)$. Here $U_0^{\delta}(\gamma_0,t)=\Bigl(u_0^{1,\delta}(\gamma_0,t),u_0^{2,\delta}(\gamma_0,t)\Bigr)^{\top}$, $t>0$, is the value of the solution $U^{\delta}$ at a new boundary point of a peeled tree $\widetilde{\Omega}$.

Let us introduce the following notation:
\begin{equation*}
a(t):=u^{1,\delta}(\gamma_0,t),\quad A(t):=u^{2,\delta}(\gamma_0,t),
\end{equation*}
where $a(t)=A(t)=0$ for $t<l_1$.

Let us consider the solution $U^{f}$ of the problem (\ref{Equation1})-(\ref{Equation3}), (\ref{Equation5}) with the boundary conditions
\begin{equation*}
u^{1,f}(\gamma_1,t)=f(t),\quad u^{1,f}(\gamma_k,t)=0,\quad k=2,\cdots,L.
\end{equation*}

By Duhamel's principle, at the new boundary vertex $\gamma_0$ we obtain
\begin{equation*}
u^{1,f}(\gamma_0,t)=a(t)\ast f(t).
\end{equation*}

By definitions of the response matrices $\mathbf{R}$ and $\widetilde{\mathbf{R}}$, the following relations hold:
\begin{equation}\label{r-entr-1}
\int\limits_0^tR_{1k}(s)f(t-s)\,ds=\int\limits_0^t\widetilde{R}_{0k}(s)[a\ast f](t-s)\,ds,\quad k=L^{'}+1,\cdots,L-1.
\end{equation}

We know (e.g., see \cite{2008_Avdonin_Kurasov}) that the response entries $R_{1k}(s)$, $\widetilde{R}_{0k}(s)$, $k=L^{'}+1,\cdots,L-1$, and $a$ admit the following representations:
\begin{equation}\label{represent-1}
R_{1k}(s)=\sum_{n=1}^{\infty}r_n\delta(s-\beta_n)+\chi_{1k}(s),\quad \chi_{1k}|_{s\in(0,\beta_1)}=0,
\end{equation}
\begin{equation}\label{represent-2}
a(s)=\sum_{n=1}^{\infty}\xi_n\delta(s-\alpha_n)+\hat{a}(s),\quad \hat{a}|_{s\in(0,\alpha_1)}=0,\quad \alpha_1=l_1,\quad \xi_1\neq0,
\end{equation}
\begin{equation}\label{represent-3}
\widetilde{R}_{0k}(s)=\sum_{n=1}^{\infty}\tilde{r}_n\delta(s-\zeta_n)+\tilde{r}_{0k}(s),\quad \tilde{r}_{0k}|_{s\in(0,\zeta_1)}=0.
\end{equation}

In the above representations, the piece-wise continuous functions $\chi_{1k}(s)$, $\hat{a}(s)$ and the numbers $r_n$, $\beta_n$, $\xi_n$, $\alpha_n$ are known. Also, the sequences $\{\beta_n\}_{n=1}^{\infty}$, $\{\alpha_n\}_{n=1}^{\infty}$, and $\{\zeta_n\}_{n=1}^{\infty}$ are strictly increasing. 

We want to recover the unknown piece-wise continuous function $\tilde{r}_{0k}(s)$ and the numbers $\tilde{r}_n$, $\zeta_n$. For this purpose, we substitute expressions (\ref{represent-1}), (\ref{represent-2}), and (\ref{represent-3}) in (\ref{r-entr-1}). As a result we obtain:
\begin{equation}\label{substitution}
\begin{split}
\sum_{n=1}^{\infty}\sum_{m=1}^{\infty}\xi_n\tilde{r}_m\delta(t-\zeta_m-\alpha_n)+\sum_{n=1}^{\infty}\xi_n\tilde{r}_{0k}(t-\alpha_n)+\sum_{n=1}^{\infty}\tilde{r}_n\hat{a}(t-\zeta_n)\\
+\hat{a}(t)\ast\tilde{r}_{0k}(t)=\sum_{n=1}^{\infty}r_n\delta(t-\beta_n)+\chi_{1k}(t),\quad k=L^{'}+1,\cdots,L-1.
\end{split}
\end{equation}

By equating singular parts in equation (\ref{substitution}), we obtain
\begin{equation}\label{singular}
\sum_{n=1}^{\infty}r_n\delta(t-\beta_n)=\sum_{n=1}^{\infty}\sum_{m=1}^{\infty}\xi_n\tilde{r}_m\delta(t-\zeta_m-\alpha_n).
\end{equation}

For physical systems the impulse response $\delta(t)$ is zero for time $t<0$, and future components of the input do not contribute to the sum. Thus, let $N_{max}$ be the largest integer such that $\beta_{N_{max}}\leq t$ and let $M$ be an integer greater than $N_{max}$. Since $\beta_n\leq \zeta_m+\alpha_n$ and $\beta_n\leq \zeta_n+\alpha_m$ for every $n\geq 1$ and $m\geq 1$ (e.g., see \cite{2020_Avdonin_Zhao}), we obtain that the impulse response is not equal to zero for $t<\beta_M\leq \zeta_m+\alpha_M$ and $t<\beta_M\leq \zeta_M+\alpha_m$ for any $m\geq 1$. Hence, for any $m\geq 1$ and $M>N_{max}$, we have
\begin{equation*}
\delta(t-\zeta_m-\alpha_M)=\delta(t-\zeta_M-\alpha_m)=\delta(t-\beta_M)=0.
\end{equation*}

Therefore, the equation (\ref{singular}) can be rewritten in the form
\begin{equation}\label{finite-sum}
\sum_{n=1}^{N_{max}}r_n\delta(t-\beta_n)=\sum_{n=1}^{N_{max}}\sum_{m=1}^{N_{max}}\xi_n\tilde{r}_m\delta(t-\zeta_m-\alpha_n).
\end{equation}

Now let us show that the equation (\ref{finite-sum}) admits the following matrix representation:
\begin{equation}\label{matrix-form}
A\tilde{r}=r,
\end{equation}
where $\tilde{r}=(\tilde{r}_1,\cdots,\tilde{r}_{N_{max}})^{\top}$, $r=(r_1,\cdots,r_{N_{max}})^{\top}$, and the matrix $A$ is lower triangular with nonzero diagonal elements. Indeed, let us consider the following cases.

\textit{Case 1.} If $\beta_1=\zeta_1+\alpha_1$, then $r_1=\xi_1\tilde{r}_1$. From this it follows that $\tilde{r}_1=r_1/\xi_1$, where $\xi_1\neq0$ and we set $A_{11}=\xi_1\neq0$.

\textit{Case 2.} If $\beta_2\neq \zeta_1+\alpha_2$, then $\beta_2=\zeta_2+\alpha_1$. So, $r_2=\xi_1\tilde{r}_2$ which implies $\tilde{r}_2=r_2/\xi_1$. Based on this, we set $A_{22}=\xi_1\neq0$.

If $\beta_2=\zeta_1+\alpha_2$, but $r_2\neq \xi_1\tilde{r}_2$, then $\zeta_1+\alpha_2=\zeta_2+\alpha_1=\beta_2$. Thus, $r_2=\xi_2\tilde{r}_1+\xi_1\tilde{r}_2$. From this it follows that $\tilde{r}_2=(r_2-\xi_2\tilde{r}_1)/\xi_1$, where $\xi_1\neq0$. We set $A_{21}=\xi_2$, where $\xi_2=(r_2-\xi_1\tilde{r}_2)/\xi_1\neq 0$. If $\beta_2=\zeta_1+\alpha_1$ and $r_2=\xi_1\tilde{r}_2$, then we need to compare $\beta_3$ with $\zeta_3+\alpha_1$.

Repeating this procedure a sufficient number of times, we can generalize the definition of the matrix $A$. For every $n,m=1,\cdots,N_{max}$ such that $n>m$ we set $A_{nm}=\xi_k$ if $\beta_n=\zeta_m+\alpha_k$ for some $k$, otherwise $A_{nm}=0$. Since for every $n\geq1$ we have that $\beta_n=\zeta_n+\alpha_1$, we obtain that $A_{nn}=\xi_1\neq0$.
Now let us show that for all $m>n$ we get that $A_{nm}=0$. Indeed, for all $m>n$ we have that $\zeta_m>\zeta_n$. Hence, the equality $\beta_n=\zeta_n+\alpha_1=\zeta_m+\alpha_k$ would mean that $\alpha_1$ needs to be greater than $\alpha_k$ for some $k\geq1$. Thus we get a contradiction. Therefore, $A_{nm}=0$ for all $m>n$.

From this it follows that the matrix $A$ is lower triangular with nonzero diagonal elements and, therefore, it is invertible. Thus one can compute unknown coefficients $\tilde{r}_m$ for any $m<N_{max}<\infty$ and numbers $\zeta_m$ correspondingly.

Thus, it remains to compute the unknown piece-wise continuous function $\tilde{r}_{0k}(s)$ for all $k=L^{'}+1,\cdots,L-1$. By equating the regular parts in (\ref{substitution}), we obtain an integral equation with respect to $\tilde{r}_{0k}(s)$, which can be solved in steps:
\begin{equation}\label{integral}
\sum_{n=1}^{\infty}\tilde{r}_n\hat{a}(t-\zeta_n)\\
+\hat{a}(t)\ast\tilde{r}_{0k}(t)=\chi_{1k}(t)-\sum_{n=1}^{\infty}\xi_n\tilde{r}_{0k}(t-\alpha_n).
\end{equation}
From this it follows that the response entries $\widetilde{R}_{0k}$ are fully recovered.

To recover $\widetilde{R}_{00}(s)$ we need to use the equation
\begin{equation}\label{convolution}
\int\limits_0^t\widetilde{R}_{00}(s)[a\ast f](t-s)\,ds=[A\ast f](t)=u^{2,f}(\gamma_0,t).
\end{equation}

Then we repeat the procedure described above. We write down the expansions for $\widetilde{R}_{00}$, $A(t)$, and $a(t)$ with singular and regular parts separated, determine the singular part, and afterwords determine the regular part from the corresponding integral equation.

To recover $\widetilde{R}_{k0}(t)$, where $k=L^{'}+1,\cdots,L-1$, we consider $U^f$, the solution of the problem (\ref{Equation1})-(\ref{Equation3}), (\ref{Equation5}) with the following boundary conditions:
\begin{equation*}
u^{1,f}(\gamma_k,t)=f(t),\quad u^{1,f}(\gamma_s,t)=0,\quad s=1,\cdots,L,\quad s\neq k.
\end{equation*}

Since we know the potential functions on edges $e_{1},\cdots,e_{N^{'}-1}$, we can recover the solution of the above problem by solving the Cauchy problem on the corresponding edges with known boundary data. Indeed, let us introduce the functions
\begin{equation*}
R_s(t)=[R_{ks}\ast f](t),\quad t\geq0.
\end{equation*}

The function $U^{f}$ solves the following Cauchy problems on edges $e_1,\cdots,e_{N^{'}-1}$:
\begin{equation}\label{cp-1}
\begin{cases}
iU_t^{f}+JU_x^{f}+QU^{f}=0,\quad x\in(0,l_s),\\
u^{1,f}(\gamma_k,t)=0,\\
u^{2,f}(\gamma_s,t)=R_s(t),\quad s=1,\cdots,N^{'}-1,\\
U^{f}(x,0)=0,\quad x\in(0,l_s).
\end{cases}
\end{equation}

Using the compatibility conditions (\ref{Equation2}), (\ref{Equation3}) at the vertex $v_1$, we calculate values of $u^{1,f}(\gamma_0,t)$ and $u^{2,f}(\gamma_0,t)$ for $t>0$ at the new boundary vertex:
\begin{equation*}
a(t):=u^{1,f}(\gamma_0,t),\quad A(t):=u^{2,f}(\gamma_0,t).
\end{equation*}

Then using the definition of the response matrix $\mathbf{\widetilde{R}}(t)$, we obtain the following relations
\begin{equation*}
\int\limits_0^t\widetilde{R}_{k0}(s)f(t-s)\,ds=A(t)-[R_{00}\ast a](t).
\end{equation*}

To complete the process of recovering response entries $\widetilde{R}_{k0}(t)$ we repeat the procedure described above. We write down expansions for $\widetilde{R}_{k0}(t)$, $A(t)$, and $a(t)$ with singular and regular parts separated. Then we determine the singular part, and next we determine the regular part from the corresponding integral equation. 

Concluding this section, we state the following corollary.
\vspace{0.1cm}

\textbf{Corollary 4.5.1.} For each $k\in\{1,\cdots,N\}$ one can recover the length of all edges, the connectivity, and the potential matrix $Q(x)$ of a finite metric tree graph $\Omega$ from the response data $\mathbf{R}(t)$ for $t\in(0,T)$ with $T>2l$, where $l=\max\limits_{1\leq k\leq L-1}\{d(\gamma_k,\gamma_L)\}$, and $\{\gamma_k\}_{k=1}^{L-1}$ is the set of leaf vertices and $\gamma_L$ is a root vertex of a tree graph $\Omega$.

\end{document}